\documentclass{amsproc}
\usepackage{amssymb}
\usepackage{amsmath}
\usepackage{eurosym}
\usepackage{amsfonts}

            \advance \topmargin by -\headheight \advance
            \topmargin by -\headsep
            \evensidemargin
            \oddsidemargin
\newtheorem{theorem}{Theorem}[section]

\newtheorem{definition}[theorem]{Definition}

\numberwithin{equation}{section}
\input{tcilatex}

\begin{document}
\title[On the Dynamics of a Higher Order Nonlinear System of Difference
Equations]{On the Dynamics of a Higher Order Nonlinear System of Difference
Equations}
\thanks{}
\author[\.{I}nci Okumu\c{s}, Y\"{u}ksel Soykan]{\.{I}nci Okumu\c{s}, Y\"{u}%
ksel Soykan}
\maketitle

\begin{center}
\textsl{Zonguldak B\"{u}lent Ecevit University, Department of Mathematics, }

\textsl{Art and Science Faculty, 67100, Zonguldak, Turkey }

\textsl{e-mail: \ inci\_okumus\_90@hotmail.com}

\textsl{yuksel\_soykan@hotmail.com}
\end{center}

\textbf{Abstract.} The aim of this paper is to investigate the dynamics of a
higher order system of rational difference equations. Our concentration is
on boundedness character, the oscillatory, the existence of unbounded
solutions and the global behavior of positive solutions for the following
system of difference equations%
\begin{equation*}
x_{n+1}=A+\frac{x_{n-m}}{z_{n}},\text{ \ }y_{n+1}=A+\frac{y_{n-m}}{z_{n}},%
\text{ \ }z_{n+1}=A+\frac{z_{n-m}}{y_{n}}\text{, \ }n=0,1,...\text{,}
\end{equation*}%
where $A$ and the initial values $x_{-i}$, $y_{-i}$, $z_{-i}$, for $%
i=0,1,...,m$, are positive real numbers.



\section{Introduction}

Difference equations appear naturally as discrete analogues and as numerical
solutions of differential and delay differential equations having many
applications in economics, population biology, computer science, probability
theory, psychology and so forth. Difference equation or discrete dynamical
system is a diverse field which impact almost every branch of pure and
applied mathematics. Recently, there has been great interest in
investigating the behavior of solutions of a system of nonlinear difference
equations and discussing the asymptotic stability of their equilibrium
points. There are many papers in which systems of difference equations have
been studied, see [\ref{Bao2015}-\ref{Zhang2007}].

In [\ref{Camouzis2004}], Camouzis and and Papaschinopoulos studied the
system of difference equations%
\begin{equation}
x_{n+1}=1+\frac{x_{n}}{y_{n-m}},\ \ y_{n+1}=1+\frac{y_{n}}{x_{n-m}},\ \
n=0,1,...,  \label{equ:bbbbbbbbbbbbbb}
\end{equation}%
with initial conditions $x_{-i}$, $y_{-i}>0$, $i=-m,-m+1,...,0$, and $m$ is
positive integer.

In [\ref{Yang2005}], Yang studied the behavior of positive solutions of the
system of difference equations%
\begin{equation}
x_{n+1}=A+\frac{y_{n-1}}{x_{n-p}y_{n-q}},\text{ \ \ }y_{n+1}=A+\frac{x_{n-1}%
}{x_{n-r}y_{n-s}},\text{ \ \ }n=1,2,...\text{,}  \label{equ:asqweasde}
\end{equation}%
where $p\geq 2$, $q\geq 2$, $r\geq 2$, $s\geq 2$, $A$ is a positive
constant, and $x_{1-\max \{p,r\}}$, $x_{2-\max \{p,r\}}$, ...,$x_{0}$, $%
y_{1-\max \{q,s\}}$, $y_{2-\max \{q,s\}}$, ...,$y_{0}$ are positive real
numbers.

In [\ref{Papaschi1998}], Papaschinopoulos and Schinas considered the system
of difference equations%
\begin{equation}
x_{n+1}=A+\frac{y_{n}}{x_{n-p}},\text{ \ \ }y_{n+1}=A+\frac{x_{n}}{y_{n-q}},%
\text{ \ \ }n=0,1,...\text{,}  \label{equ:mncsrt}
\end{equation}%
where $A\in \left( 0,\infty \right) $, $p$, $q$ are positive integers and $%
x_{-p},...,x_{0}$, $y_{-q},...,y_{0}$ are positive numbers.

In [\ref{Okumus2018}], Okumu\c{s} and Soykan studied the boundedness,
persistence and periodicity of the positive solutions and the global
asymptotic stability of the positive equilibrium points of system of the
difference equations%
\begin{equation*}
x_{n+1}=A+\frac{x_{n-1}}{z_{n}},\text{ \ }y_{n+1}=A+\frac{y_{n-1}}{z_{n}},%
\text{ \ }z_{n+1}=A+\frac{z_{n-1}}{y_{n}}\text{, \ }n=0,1,...\text{,}
\end{equation*}%
where $A\in \left( 0,\infty \right) $ and initial conditions $x_{i}$, $y_{i}$%
, $z_{i}\in \left( 0,\infty \right) $, $i=-1,0$.

Motivated by all above mentioned studies and in the light of this work in [%
\ref{Okumus2018}], in this paper, we investigate the global asymptotic
stability, boundedness character and oscillatory of positive solutions of
the system of difference equations%
\begin{equation}
x_{n+1}=A+\frac{x_{n-m}}{z_{n}},\text{ \ }y_{n+1}=A+\frac{y_{n-m}}{z_{n}},%
\text{ \ }z_{n+1}=A+\frac{z_{n-m}}{y_{n}}\text{, \ }n=0,1,...\text{,}
\label{equ:zaq}
\end{equation}%
where $A$ and the initial values $x_{-i}$, $y_{-i}$, $z_{-i}$, for $%
i=0,1,...,m$, are positive real numbers and $m$ is positive integer.

\section{Preliminaries}

We recall some basic definitions that we afterwards need in the paper.

Let us introduce the discrete dynamical system:%
\begin{eqnarray}
x_{n+1} &=&f_{1}\left(
x_{n},x_{n-1},...,x_{n-k},y_{n},y_{n-1},...,y_{n-k},z_{n},z_{n-1},...,z_{n-k}\right) 
\text{,}  \notag \\
y_{n+1} &=&f_{2}\left(
x_{n},x_{n-1},...,x_{n-k},y_{n},y_{n-1},...,y_{n-k},z_{n},z_{n-1},...,z_{n-k}\right) 
\text{,}  \label{equ:lkjhfgdsa} \\
z_{n+1} &=&f_{3}\left(
x_{n},x_{n-1},...,x_{n-k},y_{n},y_{n-1},...,y_{n-k},z_{n},z_{n-1},...,z_{n-k}\right) 
\text{,}  \notag
\end{eqnarray}%
$n\in 
\mathbb{N}
$, where $f_{1}:I_{1}^{k+1}\times I_{2}^{k+1}\times I_{3}^{k+1}\rightarrow
I_{1}$, $f_{2}:I_{1}^{k+1}\times I_{2}^{k+1}\times I_{3}^{k+1}\rightarrow
I_{2}$\ and $f_{3}:I_{1}^{k+1}\times I_{2}^{k+1}\times
I_{3}^{k+1}\rightarrow I_{3}$ are continuously differentiable functions and $%
I_{1}$, $I_{2}$, $I_{3}$ are some intervals of real numbers. Also, a
solution $\{\left( x_{n},y_{n},z_{n}\right) \}_{n=-k}^{\infty }$ of system (%
\ref{equ:lkjhfgdsa}) is uniquely determined by initial values $\left(
x_{-i},y_{-i},z_{-i}\right) \in I_{1}\times I_{2}\times I_{3}$ for $i\in
\left\{ 0,1,...,k\right\} $.

\begin{definition}
An equilibrium point of system (\ref{equ:lkjhfgdsa}) is a point $\left( 
\overline{x},\overline{y},\overline{z}\right) $ that satisfies%
\begin{eqnarray*}
\overline{x} &=&f_{1}\left( \overline{x},\overline{x},...,\overline{x},%
\overline{y},\overline{y},...,\overline{y},\overline{z},\overline{z},...,%
\overline{z}\right) \text{,} \\
\overline{y} &=&f_{2}\left( \overline{x},\overline{x},...,\overline{x},%
\overline{y},\overline{y},...,\overline{y},\overline{z},\overline{z},...,%
\overline{z}\right) \text{,} \\
\overline{z} &=&f_{3}\left( \overline{x},\overline{x},...,\overline{x},%
\overline{y},\overline{y},...,\overline{y},\overline{z},\overline{z},...,%
\overline{z}\right) \text{.}
\end{eqnarray*}
\end{definition}

Together with system (\ref{equ:lkjhfgdsa}), if we consider the associated
vector map%
\begin{equation*}
F=\left(
f_{1},x_{n},x_{n-1},...,x_{n-k},f_{2},y_{n},y_{n-1},...,y_{n-k},f_{3},z_{n-1},...,z_{n-k}\right) 
\text{,}
\end{equation*}%
then the point $\left( \overline{x},\overline{y},\overline{z}\right) $ is
also called a fixed point of the vector map $F$.

\begin{definition}
Let $\left( \overline{x},\overline{y},\overline{z}\right) $ be an
equilibrium point of system (\ref{equ:lkjhfgdsa}).

\begin{description}
\item[(a)] An equilibrium point $\left( \overline{x},\overline{y},\overline{z%
}\right) $ is called stable if, for every $\varepsilon >0$; there exists $%
\delta >0$ such that for every initial value $\left(
x_{-i},y_{-i},z_{-i}\right) \in I_{1}\times I_{2}\times I_{3}$, with%
\begin{equation*}
\dsum\nolimits_{i=-k}^{0}\left\vert x_{i}-\overline{x}\right\vert <\delta ,%
\text{ }\dsum\nolimits_{i=-k}^{0}\left\vert y_{i}-\overline{y}\right\vert
<\delta \text{, }\dsum\nolimits_{i=-k}^{0}\left\vert z_{i}-\overline{z}%
\right\vert <\delta
\end{equation*}%
implying $\left\vert x_{n}-\overline{x}\right\vert <\varepsilon $, $%
\left\vert y_{n}-\overline{y}\right\vert <\varepsilon $, $\left\vert z_{n}-%
\overline{z}\right\vert <\varepsilon $for $n\in 
\mathbb{N}
$.

\item[(b)] If an equilibrium point $\left( \overline{x},\overline{y},%
\overline{z}\right) $ of system (\ref{equ:lkjhfgdsa}) is called unstable if
it is not stable.

\item[(c)] An equilibrium point $\left( \overline{x},\overline{y},\overline{z%
}\right) $ of system (\ref{equ:lkjhfgdsa})\ is called locally asymptotically
stable if, it is stable, and if in addition there exists $\gamma >0$ such
that%
\begin{equation*}
\dsum\nolimits_{i=-k}^{0}\left\vert x_{i}-\overline{x}\right\vert <\gamma ,%
\text{ }\dsum\nolimits_{i=-k}^{0}\left\vert y_{i}-\overline{y}\right\vert
<\gamma \text{, }\dsum\nolimits_{i=-k}^{0}\left\vert z_{i}-\overline{z}%
\right\vert <\gamma
\end{equation*}%
and $\left( x_{n},y_{n},z_{n}\right) \rightarrow \left( \overline{x},%
\overline{y},\overline{z}\right) $ as $n\rightarrow \infty $.

\item[(d)] An equilibrium point $\left( \overline{x},\overline{y},\overline{z%
}\right) $ of system (\ref{equ:lkjhfgdsa})\ is called a global attractor if, 
$\left( x_{n},y_{n},z_{n}\right) \rightarrow \left( \overline{x},\overline{y}%
,\overline{z}\right) $ as $n\rightarrow \infty $.

\item[(e)] An equilibrium point $\left( \overline{x},\overline{y},\overline{z%
}\right) $ of system (\ref{equ:lkjhfgdsa})\ is called globally
asymptotically stable if it is stable, and a global attractor.
\end{description}
\end{definition}

\begin{definition}
Let $\left( \overline{x},\overline{y},\overline{z}\right) $ be an
equilibrium point of the map $F$ where $f_{1}$, $f_{2}$ and $f_{3}$ are
continuously differentiable functions at $\left( \overline{x},\overline{y},%
\overline{z}\right) $. The linearized system of system (\ref{equ:lkjhfgdsa})
about the equilibrium point $\left( \overline{x},\overline{y},\overline{z}%
\right) $ is%
\begin{equation*}
X_{n+1}=F\left( X_{n}\right) =BX_{n}\text{,}
\end{equation*}%
where%
\begin{equation*}
X_{n}=\left( 
\begin{array}{c}
x_{n} \\ 
\vdots \\ 
x_{n-k} \\ 
y_{n} \\ 
\vdots \\ 
y_{n-k} \\ 
z_{n} \\ 
\vdots \\ 
z_{n-k}%
\end{array}%
\right)
\end{equation*}%
and $B$ is a Jacobian matrix of system (\ref{equ:lkjhfgdsa}) about the
equilibrium point $\left( \overline{x},\overline{y},\overline{z}\right) $.
\end{definition}

\begin{definition}
Assume that%
\begin{equation*}
X_{n+1}=F\left( X_{n}\right) ,n=0,1,...\text{,}
\end{equation*}%
be a system of difference equations such that $\overline{X}$ is a fixed
point of $F$. If no eigenvalues of the Jacobian matrix $B$ about $\overline{X%
}$ have absolute value equal to one, then $\overline{X}$ is called
hyperbolic. If there exists an eigenvalue of the Jacobian matrix $B$ about $%
\overline{X}$ with absolute value equal to one, then $\overline{X}$ is
called nonhyperbolic.
\end{definition}

\begin{theorem}[The Linearized Stability Theorem]

Assume that%
\begin{equation*}
X_{n+1}=F\left( X_{n}\right) ,n=0,1,...\text{,}
\end{equation*}%
be a system of difference equations such that $\overline{X}$ is a fixed
point of $F$.

\begin{description}
\item[(a)] If all eigenvalues of the Jacobian matrix $B$ about $\overline{X}$
lie inside the open unit disk $\left\vert \lambda \right\vert <1$, that is,
if all of them have absolute value less than one, then $\overline{X}$\ is
locally asymptotically stable.

\item[(b)] If at least one of them has a modulus greater than one, then $%
\overline{X}$\ is unstable.
\end{description}
\end{theorem}

\section{Main Results}

In this section, we prove our main results. We deal with the following cases
of $0<A<1$, $A=1$, and $A>1$.

\begin{theorem}
If $\left( \overline{x},\overline{y},\overline{z}\right) $ is a positive
equilibrium point of system (\ref{equ:zaq}), then%
\begin{equation*}
\left( \overline{x},\overline{y},\overline{z}\right) =\left\{ 
\begin{array}{cc}
\left( A+1,A+1,A+1\right) \text{,} & \text{if }A\neq 1\text{,} \\ 
\left( \mu ,\mu ,\frac{\mu }{\mu -1}\right) \text{, }\mu \in \left( 1,\infty
\right) & \text{if }A=1\text{.}%
\end{array}%
\right.
\end{equation*}
\end{theorem}

\textbf{Proof. }It is easily seen from the definition of equilibrium point
that the equilibrium points of system (\ref{equ:zaq}) are the nonnegative
solution of the equations%
\begin{equation*}
\overline{x}=A+\frac{\overline{x}}{\overline{z}}\text{, }\overline{y}=A+%
\frac{\overline{y}}{\overline{z}}\text{, }\overline{z}=A+\frac{\overline{z}}{%
\overline{y}}\text{.}
\end{equation*}%
From this, we get%
\begin{eqnarray*}
\overline{x}\overline{z} &=&A\overline{z}+\overline{x}\text{, \ }\overline{y}%
\overline{z}=A\overline{z}+\overline{y}\text{, \ }\overline{z}\overline{y}=A%
\overline{y}+\overline{z} \\
&\Rightarrow &\overline{x}\overline{z}-\overline{x}=\overline{y}\overline{z}-%
\overline{y}\text{, \ }A\overline{z}+\overline{y}=A\overline{y}+\overline{z}
\\
&\Rightarrow &\overline{x}\left( \overline{z}-1\right) =\overline{y}\left( 
\overline{z}-1\right) \text{, \ }\overline{z}\left( A-1\right) =\overline{y}%
\left( A-1\right) \text{.}
\end{eqnarray*}%
From which it follows that if $A\neq 1$,%
\begin{equation*}
\overline{x}=\overline{y}=\overline{z}=A+1\Rightarrow \left( \overline{x},%
\overline{y},\overline{z}\right) =\left( A+1,A+1,A+1\right) \text{.}
\end{equation*}%
Also, we have%
\begin{eqnarray*}
\frac{\overline{x}\overline{z}-\overline{x}}{\overline{z}} &=&A\text{, \ }%
\frac{\overline{y}\overline{z}-\overline{y}}{\overline{z}}=A\text{, \ }\frac{%
\overline{z}\overline{y}-\overline{z}}{\overline{y}}=A \\
&\Rightarrow &\frac{\overline{x}\overline{z}-\overline{x}}{\overline{z}}=%
\frac{\overline{y}\overline{z}-\overline{y}}{\overline{z}}\text{, \ }\frac{%
\overline{y}\overline{z}-\overline{y}}{\overline{z}}=\frac{\overline{z}%
\overline{y}-\overline{z}}{\overline{y}} \\
&\Rightarrow &\overline{x}\overline{z}-\overline{x}=\overline{y}\overline{z}-%
\overline{y}\text{, \ }\overline{y}^{2}\overline{z}-\overline{y}^{2}=%
\overline{z}^{2}\overline{y}-\overline{z}^{2} \\
&\Rightarrow &\overline{x}\left( \overline{z}-1\right) =\overline{y}\left( 
\overline{z}-1\right) \text{, \ }\overline{y}\overline{z}\left( \overline{y}-%
\overline{z}\right) =\left( \overline{y}-\overline{z}\right) \left( 
\overline{y}+\overline{z}\right) \text{.}
\end{eqnarray*}%
From which it follows that if $A=1$,%
\begin{equation*}
\overline{x}=\overline{y}\text{ and\ }\overline{y}\overline{z}=\overline{y}+%
\overline{z}\Rightarrow \left( \overline{x},\overline{y},\overline{z}\right)
=\left( \mu ,\mu ,\frac{\mu }{\mu -1}\right) \text{, }\mu \in \left(
1,\infty \right) \text{.}
\end{equation*}%
In that case, we have a continuous of positive equilibriums which lie on the
hyperboloid $\overline{y}\overline{z}=\overline{y}+\overline{z}$.

\begin{theorem}
Assume that $0<A<1$. Let $\left\{ \left( x_{n},y_{n},z_{n}\right) \right\} $
be an arbitrary positive solution of the system (\ref{equ:zaq}). Then, the
following statements are true.

\begin{description}
\item[(i)] If $m$ is odd and $0<x_{2k-1}<1$, $0<y_{2k-1}<1$, $0<z_{2k-1}<1$, 
$x_{2k}>\frac{1}{1-A}$, $y_{2k}>\frac{1}{1-A}$, $z_{2k}>\frac{1}{1-A}$ for $%
k=\frac{1-m}{2},\frac{3-m}{2},...,0$, then%
\begin{eqnarray*}
\lim_{n\rightarrow \infty }x_{2n} &=&\infty \text{, \ }\lim_{n\rightarrow
\infty }y_{2n}=\infty \text{, \ }\lim_{n\rightarrow \infty }z_{2n}=\infty 
\text{,} \\
\lim_{n\rightarrow \infty }x_{2n+1} &=&A\text{, \ }\lim_{n\rightarrow \infty
}y_{2n+1}=A\text{, \ }\lim_{n\rightarrow \infty }z_{2n+1}=A\text{.}
\end{eqnarray*}

\item[(ii)] If $m$ is odd and $0<x_{2k}<1$, $0<y_{2k}<1$, $0<z_{2k}<1$, $%
x_{2k-1}>\frac{1}{1-A}$, $y_{2k-1}>\frac{1}{1-A}$, $z_{2k-1}>\frac{1}{1-A}$
for $k=\frac{1-m}{2},\frac{3-m}{2},...,0$, then%
\begin{eqnarray*}
\lim_{n\rightarrow \infty }x_{2n} &=&A\text{, \ }\lim_{n\rightarrow \infty
}y_{2n}=A\text{, \ }\lim_{n\rightarrow \infty }z_{2n}=A\text{,} \\
\lim_{n\rightarrow \infty }x_{2n+1} &=&\infty \text{, \ }\lim_{n\rightarrow
\infty }y_{2n+1}=\infty \text{, \ }\lim_{n\rightarrow \infty
}z_{2n+1}=\infty \text{.}
\end{eqnarray*}

\item[(iii)] If $m$ is even, we can not get some useful results.
\end{description}
\end{theorem}

\textbf{Proof.}

\begin{description}
\item[(i)] Clearly, we get%
\begin{eqnarray*}
0 &<&x_{1}=A+\frac{x_{-m}}{z_{0}}<A+\frac{1}{z_{0}}<A+\left( 1-A\right) =1%
\text{,} \\
0 &<&y_{1}=A+\frac{y_{-m}}{z_{0}}<A+\frac{1}{z_{0}}<A+\left( 1-A\right) =1%
\text{,} \\
0 &<&z_{1}=A+\frac{z_{-m}}{y_{0}}<A+\frac{1}{y_{0}}<A+\left( 1-A\right) =1%
\text{,} \\
x_{2} &=&A+\frac{x_{1-m}}{z_{1}}>x_{1-m}>\frac{1}{1-A}\text{,} \\
y_{2} &=&A+\frac{y_{1-m}}{z_{1}}>y_{1-m}>\frac{1}{1-A}\text{,} \\
z_{2} &=&A+\frac{z_{1-m}}{y_{1}}>z_{1-m}>\frac{1}{1-A}\text{.}
\end{eqnarray*}%
By induction for $n=1,2,...,$ we obtain%
\begin{eqnarray}
0 &<&x_{2n-1}<1\text{, }0<y_{2n-1}<1\text{, }0<z_{2n-1}<1\text{,}
\label{ineq:lokij} \\
x_{2n} &>&\frac{1}{1-A}\text{, }y_{2n}>\frac{1}{1-A}\text{, }z_{2n}>\frac{1}{%
1-A}\text{.}  \notag
\end{eqnarray}%
Thus, for $n\geq \left( m+2\right) /2$,%
\begin{eqnarray*}
x_{2n} &=&A+\frac{x_{2n-\left( m+1\right) }}{z_{2n-1}}>A+x_{2n-\left(
m+1\right) }=2A+\frac{x_{2n-\left( 2m+2\right) }}{z_{2n-\left( m+2\right) }}
\\
&>&2A+x_{2n-\left( 2m+2\right) }\text{,} \\
y_{2n} &=&A+\frac{y_{2n-\left( m+1\right) }}{z_{2n-1}}>A+y_{2n-\left(
m+1\right) }=2A+\frac{y_{2n-\left( 2m+2\right) }}{z_{2n-\left( m+2\right) }}
\\
&>&2A+y_{2n-\left( 2m+2\right) }\text{,} \\
z_{2n} &=&A+\frac{z_{2n-\left( m+1\right) }}{y_{2n-1}}>A+z_{2n-\left(
m+1\right) }=2A+\frac{z_{2n-\left( 2m+2\right) }}{y_{2n-\left( m+2\right) }}
\\
&>&2A+z_{2n-\left( 2m+2\right) }\text{,}
\end{eqnarray*}%
from which we get%
\begin{equation*}
\lim_{n\rightarrow \infty }x_{2n}=\infty \text{, \ }\lim_{n\rightarrow
\infty }y_{2n}=\infty \text{, \ }\lim_{n\rightarrow \infty }z_{2n}=\infty 
\text{.}
\end{equation*}%
Noting that (\ref{ineq:lokij}) and taking limits on the both sides of three
equations%
\begin{equation*}
x_{2n+1}=A+\frac{x_{2n-m}}{z_{2n}}\text{, }y_{2n+1}=A+\frac{y_{2n-m}}{z_{2n}}%
\text{, }z_{2n+1}=A+\frac{z_{2n-m}}{y_{2n}}\text{,}
\end{equation*}%
we have%
\begin{equation*}
\lim_{n\rightarrow \infty }x_{2n+1}=A\text{, \ }\lim_{n\rightarrow \infty
}y_{2n+1}=A\text{, \ }\lim_{n\rightarrow \infty }z_{2n+1}=A\text{.}
\end{equation*}

\item[(ii)] Obviously, we have%
\begin{eqnarray*}
x_{1} &=&A+\frac{x_{-m}}{z_{0}}>x_{-m}>\frac{1}{1-A}\text{,} \\
y_{1} &=&A+\frac{y_{-m}}{z_{0}}>y_{-m}>\frac{1}{1-A}\text{,} \\
z_{1} &=&A+\frac{z_{-m}}{y_{0}}>z_{-m}>\frac{1}{1-A}\text{,} \\
0 &<&x_{2}=A+\frac{x_{1-m}}{z_{1}}<A+\frac{1}{z_{1}}<A+\left( 1-A\right) =1%
\text{,} \\
0 &<&y_{2}=A+\frac{y_{1-m}}{z_{1}}<A+\frac{1}{z_{1}}<A+\left( 1-A\right) =1%
\text{,} \\
0 &<&z_{2}=A+\frac{z_{1-m}}{y_{1}}\text{.}<A+\frac{1}{y_{1}}<A+\left(
1-A\right) =1\text{.}
\end{eqnarray*}%
By induction for $n=1,2,...,$ we obtain%
\begin{eqnarray}
x_{2n-1} &>&\frac{1}{1-A}\text{, }y_{2n-1}>\frac{1}{1-A}\text{, }z_{2n-1}>%
\frac{1}{1-A}\text{,}  \label{ineq:opnmbv} \\
0 &<&x_{2n}<1,0<y_{2n}<1,0<z_{2n}<1.  \notag
\end{eqnarray}%
So, for $n\geq \left( m+2\right) /2$,%
\begin{eqnarray*}
x_{2n+1} &=&A+\frac{x_{2n-m}}{z_{2n}}>A+x_{2n-m}=2A+\frac{x_{\left(
2n-2m\right) -1}}{z_{2n-\left( m+1\right) }} \\
&>&2A+x_{\left( 2n-2m\right) -1}\text{,} \\
y_{2n+1} &=&A+\frac{y_{2n-m}}{z_{2n}}>A+y_{2n-m}=2A+\frac{y_{\left(
2n-2m\right) -1}}{z_{2n-\left( m+1\right) }} \\
&>&2A+y_{\left( 2n-2m\right) -1}\text{,} \\
z_{2n+1} &=&A+\frac{z_{2n-m}}{y_{2n}}>A+z_{2n-m}=2A+\frac{z_{\left(
2n-2m\right) -1}}{y_{2n-\left( m+1\right) }} \\
&>&2A+z_{\left( 2n-2m\right) -1}\text{,}
\end{eqnarray*}%
from which we get%
\begin{equation*}
\lim_{n\rightarrow \infty }x_{2n+1}=\infty \text{, \ }\lim_{n\rightarrow
\infty }y_{2n+1}=\infty \text{, \ }\lim_{n\rightarrow \infty
}z_{2n+1}=\infty \text{.}
\end{equation*}%
Noting that (\ref{ineq:opnmbv}) and taking limits on the both sides of three
equations%
\begin{equation*}
x_{2n}=A+\frac{x_{2n-\left( m+1\right) }}{z_{2n-1}}\text{, }y_{2n+1}=A+\frac{%
y_{2n-\left( m+1\right) }}{z_{2n-1}}\text{, }z_{2n+1}=A+\frac{z_{2n-\left(
m+1\right) }}{y_{2n-1}}\text{,}
\end{equation*}%
we have%
\begin{equation*}
\lim_{n\rightarrow \infty }x_{2n}=A\text{, \ }\lim_{n\rightarrow \infty
}y_{2n}=A\text{, \ }\lim_{n\rightarrow \infty }z_{2n+}=A\text{.}
\end{equation*}
\end{description}

\begin{theorem}
Suppose that $A=1$. Then every positive solution of system (\ref{equ:zaq})
is bounded and persists.
\end{theorem}

\textbf{Proof.} Let $\left\{ \left( x_{n},y_{n},z_{n}\right) \right\} $ be a
positive solution of the system (\ref{equ:zaq}).

Obviously, $x_{n}>1$, $y_{n}>1$, $z_{n}>1$, for $n\geq 1$. So, we have%
\begin{equation*}
x_{i},y_{i},z_{i}\in \left[ M,\frac{M}{M-1}\right] \text{, \ }i=1,2,...,m+1%
\text{,}
\end{equation*}%
where%
\begin{equation*}
M=\min \left\{ \alpha ,\frac{\beta }{\beta -1}\right\} >1,\alpha
=\min_{1\leq i\leq m+1}\left\{ x_{i},y_{i},z_{i}\right\} ,\beta =\max_{1\leq
i\leq m+1}\left\{ x_{i},y_{i},z_{i}\right\} \text{.}
\end{equation*}

Then, we obtain%
\begin{eqnarray*}
M &=&1+\frac{M}{M/\left( M-1\right) }\leq x_{m+2}=1+\frac{x_{1}}{z_{m+1}}%
\leq 1+\frac{M/\left( M-1\right) }{M}=\frac{M}{M-1}\text{,} \\
M &=&1+\frac{M}{M/\left( M-1\right) }\leq y_{m+2}=1+\frac{y_{1}}{z_{m+1}}%
\leq 1+\frac{M/\left( M-1\right) }{M}=\frac{M}{M-1}\text{,} \\
M &=&1+\frac{M}{M/\left( M-1\right) }\leq z_{m+2}=1+\frac{z_{1}}{y_{m+1}}%
\leq 1+\frac{M/\left( M-1\right) }{M}=\frac{M}{M-1}\text{.}
\end{eqnarray*}%
By induction, we get%
\begin{equation*}
x_{i},y_{i},z_{i}\in \left[ M,\frac{M}{M-1}\right] \text{, \ }i=1,2,...\text{%
.}
\end{equation*}

\begin{theorem}
Assume that $A=1$. Let $\left\{ \left( x_{n},y_{n},z_{n}\right) \right\} $
be a positive solution of the system (\ref{equ:zaq}). Then, either $\left\{
\left( x_{n},y_{n},z_{n}\right) \right\} $ consists of a single semicycle or 
$\left\{ \left( x_{n},y_{n},z_{n}\right) \right\} $ oscillates about the
equilibrium point $\left( \overline{x},\overline{y},\overline{z}\right)
=\left( \mu ,\mu ,\frac{\mu }{\mu -1}\right) $ with semicycles having at
most $m$ terms$.$
\end{theorem}

\textbf{Proof.} Suppose that $\left\{ \left( x_{n},y_{n},z_{n}\right)
\right\} $ has at least two semicycles \ Then, there exists $N\geq -m$ such
that either $x_{N}<\overline{x}\leq x_{N+1}$ or $x_{N+1}<\overline{x}\leq
x_{N}$ ($y_{N}<\overline{y}\leq y_{N+1}$ or $y_{N+1}<\overline{y}\leq y_{N}$
and $z_{N}<\overline{z}\leq z_{N+1}$ or $z_{N+1}<\overline{z}\leq z_{N}$).
Firstly, we assume that the case $x_{N}<\overline{x}\leq x_{N+1}$, $y_{N}<%
\overline{y}\leq y_{N+1}$ and $z_{N}<\overline{z}\leq z_{N+1}$. Since the
other case is similar, it will be omitted. Suppose that rhe positive
semicycle beginning with the term $\left( x_{N+1}\text{, }y_{N+1}\text{, }%
z_{N+1}\right) $ have $m$ terms. Then we have%
\begin{eqnarray*}
x_{N+1} &<&\overline{x}=\mu \leq x_{N+m}\text{,} \\
y_{N+1} &<&\overline{y}=\mu \leq y_{N+m}\text{,} \\
z_{N+1} &<&\overline{z}=\frac{\mu }{\mu -1}\leq z_{N+m}\text{.}
\end{eqnarray*}%
Therefore, we get%
\begin{eqnarray*}
x_{N+m+1} &=&1+\frac{x_{N}}{z_{N+m}}<1+\frac{\overline{x}}{\overline{z}}=\mu 
\text{,} \\
y_{N+m+1} &=&1+\frac{y_{N}}{z_{N+m}}<1+\frac{\overline{y}}{\overline{z}}=\mu 
\text{,} \\
z_{N+m+1} &=&1+\frac{z_{N}}{y_{N+m}}<1+\frac{\overline{z}}{\overline{y}}=%
\frac{\mu }{\mu -1}\text{.}
\end{eqnarray*}%
This completes the proof.

\begin{theorem}
\label{teo:jhgfdswedf}Suppose that $A>1$. Then every positive solution of
system (\ref{equ:zaq}) is bounded and persists.
\end{theorem}

\textbf{Proof.} Let $\left\{ \left( x_{n},y_{n},z_{n}\right) \right\} $ be a
positive solution of the system (\ref{equ:zaq}).

Obviously, $x_{n}>A>1$, $y_{n}>A>1$, $z_{n}>A>1$, for $n\geq 1$. So, we have%
\begin{equation*}
x_{i},y_{i},z_{i}\in \left[ M,\frac{M}{M-A}\right] \text{, \ }i=1,2,...,m+1%
\text{,}
\end{equation*}%
where%
\begin{equation*}
M=\min \left\{ \alpha ,\frac{\beta }{\beta -1}\right\} >1,\alpha
=\min_{1\leq i\leq m+1}\left\{ x_{i},y_{i},z_{i}\right\} ,\beta =\max_{1\leq
i\leq m+1}\left\{ x_{i},y_{i},z_{i}\right\} \text{.}
\end{equation*}

Then, we obtain%
\begin{eqnarray*}
M &=&A+\frac{M}{M/\left( M-A\right) }\leq x_{m+2}=1+\frac{x_{1}}{z_{m+1}}%
\leq 1+\frac{M/\left( M-A\right) }{M}=\frac{M}{M-A}\text{,} \\
M &=&A+\frac{M}{M/\left( M-A\right) }\leq y_{m+2}=1+\frac{y_{1}}{z_{m+1}}%
\leq 1+\frac{M/\left( M-A\right) }{M}=\frac{M}{M-A}\text{,} \\
M &=&A+\frac{M}{M/\left( M-A\right) }\leq z_{m+2}=1+\frac{z_{1}}{y_{m+1}}%
\leq 1+\frac{M/\left( M-A\right) }{M}=\frac{M}{M-A}\text{.}
\end{eqnarray*}%
By induction, we get%
\begin{equation*}
x_{i},y_{i},z_{i}\in \left[ M,\frac{M}{M-A}\right] \text{, \ }i=1,2,...\text{%
.}
\end{equation*}%
The proof is completed.

Before we give the following theorems about the stability of the equilibrium
points, we consider the following transformation to build the corresponding
linearized form of system (\ref{equ:zaq}) :%
\begin{eqnarray*}
&&\left(
x_{n},x_{n-1},...,x_{n-m},y_{n},y_{n-1},...,y_{n-m},z_{n},z_{n-1},...,z_{n-m}\right)
\\
&\rightarrow &\left(
f,f_{1},...,f_{m},g,g_{1},...,g_{m},h,h_{1},...,h_{m}\right)
\end{eqnarray*}%
where%
\begin{eqnarray*}
f &=&A+\frac{x_{n-m}}{z_{n}} \\
f_{1} &=&x_{n} \\
&&\vdots \\
f_{m} &=&x_{n-m} \\
g &=&A+\frac{y_{n-m}}{z_{n}} \\
g_{1} &=&y_{n} \\
&&\vdots \\
g_{m} &=&y_{n-m} \\
h &=&A+\frac{z_{n-m}}{y_{n}} \\
h_{1} &=&z_{n} \\
&&\vdots \\
h_{m} &=&z_{n-m}\text{.}
\end{eqnarray*}%
The Jacobian matrix about the equilibrium point $\left( \overline{x},%
\overline{y},\overline{z}\right) $ under the above transformation is given by%
\begin{equation*}
B\left( \overline{x},\overline{y},\overline{z}\right) =\left( 
\begin{array}{cccccccccccc}
0 & \ldots & 0 & \frac{1}{\overline{z}} & 0 & \ldots & 0 & 0 & -\frac{%
\overline{x}}{\overline{z}^{2}} & \ldots & 0 & 0 \\ 
1 & \ldots & 0 & 0 & 0 & \ldots & 0 & 0 & 0 & \ldots & 0 & 0 \\ 
\vdots & \ddots & \vdots & \vdots & \vdots & \ddots & \vdots & \vdots & 
\vdots & \ddots & \vdots & \vdots \\ 
0 & \ldots & 1 & 0 & 0 & \ldots & 0 & 0 & 0 & \ldots & 0 & 0 \\ 
0 & \ldots & 0 & 0 & 0 & \ldots & 0 & \frac{1}{\overline{z}} & -\frac{%
\overline{y}}{\overline{z}^{2}} & \ldots & 0 & 0 \\ 
0 & \ldots & 0 & 0 & 1 & \ldots & 0 & 0 & 0 & \ldots & 0 & 0 \\ 
\vdots & \ddots & \vdots & \vdots & \vdots & \ddots & \vdots & \vdots & 
\vdots & \ddots & \vdots & \vdots \\ 
0 & \ldots & 0 & 0 & 0 & \ldots & 1 & 0 & 0 & \ldots & 0 & 0 \\ 
0 & \ldots & 0 & 0 & -\frac{\overline{z}}{\overline{y}^{2}} & \ldots & 0 & 0
& 0 & \ldots & 0 & \frac{1}{\overline{y}} \\ 
0 & \ldots & 0 & 0 & 0 & \ldots & 0 & 0 & 1 & \ldots & 0 & 0 \\ 
\vdots & \ddots & \vdots & \vdots & \vdots & \ddots & \vdots & \vdots & 
\vdots & \ddots & \vdots & \vdots \\ 
0 & \ldots & 0 & 0 & 0 & \ldots & 0 & 0 & 0 & \ldots & 1 & 0%
\end{array}%
\right) \text{,}
\end{equation*}%
where $B=\left( b_{ij}\right) $, $1\leq i,j\leq 3m+3$ is an $\left(
3m+3\right) \times \left( 3m+3\right) $ matrix.

\begin{theorem}
If $A=1$, then the equilibrium point of system (\ref{equ:zaq}) is locally
asymptotically stable.
\end{theorem}

\textbf{Proof. }The linearized system of system (\ref{equ:zaq}) about the
equilibrium point $\left( \overline{x},\overline{y},\overline{z}\right)
=\left( \mu ,\mu ,\frac{\mu }{\mu -1}\right) $ is%
\begin{equation*}
X_{n+1}=BX_{n}\text{,}
\end{equation*}%
where $X_{n}=\left(
x_{n},x_{n-1},...,x_{n-m},y_{n},y_{n-1},...,y_{n-m},z_{n},z_{n-1},...,z_{n-m}\right) ^{T} 
$ and%
\begin{equation*}
B\left( \overline{x},\overline{y},\overline{z}\right) =\left( 
\begin{array}{cccccccccccc}
0 & \ldots & 0 & \frac{\mu -1}{\mu } & 0 & \ldots & 0 & 0 & -\frac{\left(
\mu -1\right) ^{2}}{\mu } & \ldots & 0 & 0 \\ 
1 & \ldots & 0 & 0 & 0 & \ldots & 0 & 0 & 0 & \ldots & 0 & 0 \\ 
\vdots & \ddots & \vdots & \vdots & \vdots & \ddots & \vdots & \vdots & 
\vdots & \ddots & \vdots & \vdots \\ 
0 & \ldots & 1 & 0 & 0 & \ldots & 0 & 0 & 0 & \ldots & 0 & 0 \\ 
0 & \ldots & 0 & 0 & 0 & \ldots & 0 & \frac{\mu -1}{\mu } & -\frac{\left(
\mu -1\right) ^{2}}{\mu } & \ldots & 0 & 0 \\ 
0 & \ldots & 0 & 0 & 1 & \ldots & 0 & 0 & 0 & \ldots & 0 & 0 \\ 
\vdots & \ddots & \vdots & \vdots & \vdots & \ddots & \vdots & \vdots & 
\vdots & \ddots & \vdots & \vdots \\ 
0 & \ldots & 0 & 0 & 0 & \ldots & 1 & 0 & 0 & \ldots & 0 & 0 \\ 
0 & \ldots & 0 & 0 & -\frac{1}{\mu \left( \mu -1\right) ^{2}} & \ldots & 0 & 
0 & 0 & \ldots & 0 & \frac{1}{\mu } \\ 
0 & \ldots & 0 & 0 & 0 & \ldots & 0 & 0 & 1 & \ldots & 0 & 0 \\ 
\vdots & \ddots & \vdots & \vdots & \vdots & \ddots & \vdots & \vdots & 
\vdots & \ddots & \vdots & \vdots \\ 
0 & \ldots & 0 & 0 & 0 & \ldots & 0 & 0 & 0 & \ldots & 1 & 0%
\end{array}%
\right) \text{.}
\end{equation*}

Let $\lambda _{1},\lambda _{2},...,\lambda _{3m+3}$ denote the $3m+3$
eigenvalues of the matrix $B$ and $D=diag\left(
d_{1},d_{2},...,d_{3m+3}\right) $ be a diagonal matrix, where $%
d_{1}=d_{m+2}=d_{2m+3}=1$, $d_{1+k}=d_{m+2+k}=d_{2m+3+k}=1-k\varepsilon $, $%
1\leq k\leq m$ and%
\begin{equation*}
0<\varepsilon <\left\{ \frac{\mu ^{2}-2\mu +2}{m\mu },\frac{\mu ^{2}-2\mu +2%
}{m\mu \left( \mu -1\right) }\right\} \text{.}
\end{equation*}%
Obviously, $D$ is invertible. Computing matrix $DBD^{-1}$, we have that%
\begin{equation*}
DBD^{-1}
\end{equation*}%
is equal to

$\left( 
\begin{array}{cccccccccccc}
0 & \ldots  & 0 & \frac{\mu -1}{\mu }\frac{d_{1}}{d_{m+1}} & 0 & \ldots  & 0
& 0 & -\frac{\left( \mu -1\right) ^{2}}{\mu }\frac{d_{1}}{d_{2m+3}} & \ldots 
& 0 & 0 \\ 
\frac{d_{2}}{d_{1}} & \ldots  & 0 & 0 & 0 & \ldots  & 0 & 0 & 0 & \ldots  & 0
& 0 \\ 
\vdots  & \ddots  & \vdots  & \vdots  & \vdots  & \ddots  & \vdots  & \vdots 
& \vdots  & \ddots  & \vdots  & \vdots  \\ 
0 & \ldots  & \frac{d_{m+1}}{d_{m}} & 0 & 0 & \ldots  & 0 & 0 & 0 & \ldots 
& 0 & 0 \\ 
0 & \ldots  & 0 & 0 & 0 & \ldots  & 0 & \frac{\mu -1}{\mu }\frac{d_{m+2}}{%
d_{2m+2}} & -\frac{\left( \mu -1\right) ^{2}}{\mu }\frac{d_{m+2}}{d_{2m+3}}
& \ldots  & 0 & 0 \\ 
0 & \ldots  & 0 & 0 & \frac{d_{m+3}}{d_{m+2}} & \ldots  & 0 & 0 & 0 & \ldots 
& 0 & 0 \\ 
\vdots  & \ddots  & \vdots  & \vdots  & \vdots  & \ddots  & \vdots  & \vdots 
& \vdots  & \ddots  & \vdots  & \vdots  \\ 
0 & \ldots  & 0 & 0 & 0 & \ldots  & \frac{d_{2m+2}}{d_{2m+1}} & 0 & 0 & 
\ldots  & 0 & 0 \\ 
0 & \ldots  & 0 & 0 & -\frac{1}{\mu \left( \mu -1\right) ^{2}}\frac{d_{2m+3}%
}{d_{m+2}} & \ldots  & 0 & 0 & 0 & \ldots  & 0 & \frac{1}{\mu }\frac{d_{2m+3}%
}{d_{3m+3}} \\ 
0 & \ldots  & 0 & 0 & 0 & \ldots  & 0 & 0 & \frac{d_{2m+4}}{d_{2m+3}} & 
\ldots  & 0 & 0 \\ 
\vdots  & \ddots  & \vdots  & \vdots  & \vdots  & \ddots  & \vdots  & \vdots 
& \vdots  & \ddots  & \vdots  & \vdots  \\ 
0 & \ldots  & 0 & 0 & 0 & \ldots  & 0 & 0 & 0 & \ldots  & \frac{d_{3m+3}}{%
d_{3m+2}} & 0%
\end{array}%
\right) $.

The three chains of inequalities%
\begin{eqnarray*}
1 &=&d_{1}>d_{2}>...>d_{m}>d_{m+1}>0\text{,} \\
1 &=&d_{m+2}>d_{m+3}>...>d_{2m+1}>d_{2m+2}>0\text{,} \\
1 &=&d_{2m+3}>d_{2m+4}>...>d_{3m+2}>d_{3m+3}>0\text{,}
\end{eqnarray*}%
imply that%
\begin{eqnarray*}
d_{2}d_{1}^{-1} &<&1\text{, }d_{3}d_{2}^{-1}<1\text{,}...\text{,}%
d_{m+1}d_{m}^{-1}<1\text{,} \\
d_{m+3}d_{m+2}^{-1} &<&1\text{, }d_{m+4}d_{m+3}^{-1}<1\text{,}...\text{,}%
d_{2m+2}d_{2m+1}^{-1}<1\text{,} \\
d_{2m+4}d_{2m+3}^{-1} &<&1\text{, }d_{2m+5}d_{2m+4}^{-1}<1\text{,}...\text{,}%
d_{3m+3}d_{3m+2}^{-1}<1\text{.}
\end{eqnarray*}%
Also,%
\begin{eqnarray*}
\left( \frac{\mu -1}{\mu }\right) d_{1}d_{m+1}^{-1}+\left( -\frac{\left( \mu
-1\right) ^{2}}{\mu }\right) d_{1}d_{2m+3}^{-1} &=&\left( \frac{\mu -1}{\mu }%
\right) \left( \frac{1}{1-m\varepsilon }\right) -\frac{\left( \mu -1\right)
^{2}}{\mu } \\
&<&\left( \frac{\mu -1}{\mu }\right) \left( \frac{1}{1-m\varepsilon }\right)
-\left( \frac{\left( \mu -1\right) ^{2}}{\mu }\right) \left( \frac{1}{%
1-m\varepsilon }\right) \\
&=&\left( \frac{-\mu ^{2}+3\mu -2}{\mu }\right) \left( \frac{1}{%
1-m\varepsilon }\right) <1\text{,}
\end{eqnarray*}%
\begin{eqnarray*}
\left( \frac{\mu -1}{\mu }\right) d_{m+2}d_{2m+2}^{-1}+\left( -\frac{\left(
\mu -1\right) ^{2}}{\mu }\right) d_{m+2}d_{2m+3}^{-1} &=&\left( \frac{\mu -1%
}{\mu }\right) \left( \frac{1}{1-m\varepsilon }\right) -\frac{\left( \mu
-1\right) ^{2}}{\mu } \\
&=&\left( \frac{-\mu ^{2}+3\mu -2}{\mu }\right) \left( \frac{1}{%
1-m\varepsilon }\right) <1\text{,}
\end{eqnarray*}%
\begin{eqnarray*}
\left( -\frac{1}{\mu \left( \mu -1\right) ^{2}}\right)
d_{2m+3}d_{m+2}^{-1}+\left( \frac{1}{\mu }\right) d_{2m+3}d_{3m+3}^{-1}
&=&\left( -\frac{1}{\mu \left( \mu -1\right) ^{2}}\right) +\left( \frac{1}{%
\mu }\right) \left( \frac{1}{1-m\varepsilon }\right) \\
&<&\left( -\frac{1}{\mu \left( \mu -1\right) ^{2}}\right) \left( \frac{1}{%
1-m\varepsilon }\right) +\left( \frac{1}{\mu }\right) \left( \frac{1}{%
1-m\varepsilon }\right) \\
&=&\left( \frac{\mu -2}{\mu \left( \mu -1\right) }\right) \left( \frac{1}{%
1-m\varepsilon }\right) <1\text{.}
\end{eqnarray*}%
Since $B$ has the same eigenvalues as $DBD^{-1}=E=\left( e_{ij}\right) $, we
obtain that%
\begin{eqnarray*}
\max_{1\leq i\leq 3m+3}\left\vert \lambda _{i}\right\vert &\leq &\left\Vert
DBD^{-1}\right\Vert _{\infty } \\
&=&\max_{1\leq i\leq 3m+3}\left\{ \dsum\limits_{j=1}^{3m+3}\left\vert
e_{ij}\right\vert \right\} \\
&=&\max \left\{ 
\begin{array}{c}
d_{2}d_{1}^{-1}\text{, }d_{3}d_{2}^{-1}\text{,}...\text{,}d_{m+1}d_{m}^{-1}%
\text{,} \\ 
d_{m+3}d_{m+2}^{-1}\text{, }d_{m+4}d_{m+3}^{-1}\text{,}...\text{,}%
d_{2m+2}d_{2m+1}^{-1}\text{,} \\ 
d_{2m+4}d_{2m+3}^{-1}\text{, }d_{2m+5}d_{2m+4}^{-1}\text{,}...\text{,}%
d_{3m+3}d_{3m+2}^{-1}\text{,} \\ 
\left( \frac{\mu -1}{\mu }\right) d_{1}d_{m+1}^{-1}-\left( \frac{\left( \mu
-1\right) ^{2}}{\mu }\right) d_{1}d_{2m+3}^{-1}\text{,} \\ 
\left( \frac{\mu -1}{\mu }\right) d_{m+2}d_{2m+2}^{-1}-\left( \frac{\left(
\mu -1\right) ^{2}}{\mu }\right) d_{m+2}d_{2m+3}^{-1}\text{,} \\ 
\left( -\frac{1}{\mu \left( \mu -1\right) ^{2}}\right)
d_{2m+3}d_{m+2}^{-1}+\left( \frac{1}{\mu }\right) d_{2m+3}d_{3m+3}^{-1}%
\end{array}%
\right\} \\
&<&1\text{.}
\end{eqnarray*}%
This implies that the equilibrium point of system (\ref{equ:zaq}) is locally
asymptotically stable.

\begin{theorem}
If $A>1$, then the equilibrium point of system (\ref{equ:zaq}) is locally
asymptotically stable.
\end{theorem}

\textbf{Proof. }The linearized system of system (\ref{equ:zaq}) about the
equilibrium point $\left( \overline{x},\overline{y},\overline{z}\right)
=\left( A+1,A+1,A+1\right) $ is%
\begin{equation*}
X_{n+1}=BX_{n}\text{,}
\end{equation*}%
where $X_{n}=\left(
x_{n},x_{n-1},...,x_{n-m},y_{n},y_{n-1},...,y_{n-m},z_{n},z_{n-1},...,z_{n-m}\right) ^{T} 
$ and%
\begin{equation*}
B\left( \overline{x},\overline{y},\overline{z}\right) =\left( 
\begin{array}{cccccccccccc}
0 & \ldots & 0 & c^{-1} & 0 & \ldots & 0 & 0 & -c^{-1} & \ldots & 0 & 0 \\ 
1 & \ldots & 0 & 0 & 0 & \ldots & 0 & 0 & 0 & \ldots & 0 & 0 \\ 
\vdots & \ddots & \vdots & \vdots & \vdots & \ddots & \vdots & \vdots & 
\vdots & \ddots & \vdots & \vdots \\ 
0 & \ldots & 1 & 0 & 0 & \ldots & 0 & 0 & 0 & \ldots & 0 & 0 \\ 
0 & \ldots & 0 & 0 & 0 & \ldots & 0 & c^{-1} & -c^{-1} & \ldots & 0 & 0 \\ 
0 & \ldots & 0 & 0 & 1 & \ldots & 0 & 0 & 0 & \ldots & 0 & 0 \\ 
\vdots & \ddots & \vdots & \vdots & \vdots & \ddots & \vdots & \vdots & 
\vdots & \ddots & \vdots & \vdots \\ 
0 & \ldots & 0 & 0 & 0 & \ldots & 1 & 0 & 0 & \ldots & 0 & 0 \\ 
0 & \ldots & 0 & 0 & -c^{-1} & \ldots & 0 & 0 & 0 & \ldots & 0 & c^{-1} \\ 
0 & \ldots & 0 & 0 & 0 & \ldots & 0 & 0 & 1 & \ldots & 0 & 0 \\ 
\vdots & \ddots & \vdots & \vdots & \vdots & \ddots & \vdots & \vdots & 
\vdots & \ddots & \vdots & \vdots \\ 
0 & \ldots & 0 & 0 & 0 & \ldots & 0 & 0 & 0 & \ldots & 1 & 0%
\end{array}%
\right) \text{,}
\end{equation*}%
where $c=A+1$.

Let $\lambda _{1},\lambda _{2},...,\lambda _{3m+3}$ denote the $3m+3$
eigenvalues of the matrix $B$ and $D=diag\left(
d_{1},d_{2},...,d_{3m+3}\right) $ be a diagonal matrix, where $%
d_{1}=d_{m+2}=d_{2m+3}=1$, $d_{1+k}=d_{m+2+k}=d_{2m+3+k}=1-k\varepsilon $, $%
1\leq k\leq m$ and%
\begin{equation*}
0<\varepsilon <\left\{ \frac{1}{m},\frac{c-2}{cm}\right\} \text{.}
\end{equation*}%
Obviously, $D$ is invertible. Computing matrix $DBD^{-1}$, we have that%
\begin{equation*}
DBD^{-1}
\end{equation*}

$=\left( 
\begin{array}{cccccccccccc}
0 & \ldots & 0 & \frac{c^{-1}d_{1}}{d_{m+1}} & 0 & \ldots & 0 & 0 & \frac{%
-c^{-1}d_{1}}{d_{2m+3}} & \ldots & 0 & 0 \\ 
\frac{d_{2}}{d_{1}} & \ldots & 0 & 0 & 0 & \ldots & 0 & 0 & 0 & \ldots & 0 & 
0 \\ 
\vdots & \ddots & \vdots & \vdots & \vdots & \ddots & \vdots & \vdots & 
\vdots & \ddots & \vdots & \vdots \\ 
0 & \ldots & \frac{d_{m+1}}{d_{m}} & 0 & 0 & \ldots & 0 & 0 & 0 & \ldots & 0
& 0 \\ 
0 & \ldots & 0 & 0 & 0 & \ldots & 0 & \frac{c^{-1}d_{m+2}}{d_{2m+2}} & \frac{%
-c^{-1}d_{m+2}}{d_{2m+3}} & \ldots & 0 & 0 \\ 
0 & \ldots & 0 & 0 & \frac{d_{m+3}}{d_{m+2}} & \ldots & 0 & 0 & 0 & \ldots & 
0 & 0 \\ 
\vdots & \ddots & \vdots & \vdots & \vdots & \ddots & \vdots & \vdots & 
\vdots & \ddots & \vdots & \vdots \\ 
0 & \ldots & 0 & 0 & 0 & \ldots & \frac{d_{2m+2}}{d_{2m+1}} & 0 & 0 & \ldots
& 0 & 0 \\ 
0 & \ldots & 0 & 0 & \frac{-c^{-1}d_{2m+3}}{d_{m+2}} & \ldots & 0 & 0 & 0 & 
\ldots & 0 & \frac{c^{-1}d_{2m+3}}{d_{3m+3}} \\ 
0 & \ldots & 0 & 0 & 0 & \ldots & 0 & 0 & \frac{d_{2m+4}}{d_{2m+3}} & \ldots
& 0 & 0 \\ 
\vdots & \ddots & \vdots & \vdots & \vdots & \ddots & \vdots & \vdots & 
\vdots & \ddots & \vdots & \vdots \\ 
0 & \ldots & 0 & 0 & 0 & \ldots & 0 & 0 & 0 & \ldots & \frac{d_{3m+3}}{%
d_{3m+2}} & 0%
\end{array}%
\right) $.

The three chains of inequalities%
\begin{eqnarray*}
1 &=&d_{1}>d_{2}>...>d_{m}>d_{m+1}>0\text{,} \\
1 &=&d_{m+2}>d_{m+3}>...>d_{2m+1}>d_{2m+2}>0\text{,} \\
1 &=&d_{2m+3}>d_{2m+4}>...>d_{3m+2}>d_{3m+3}>0\text{,}
\end{eqnarray*}%
imply that%
\begin{eqnarray*}
d_{2}d_{1}^{-1} &<&1\text{, }d_{3}d_{2}^{-1}<1\text{,}...\text{,}%
d_{m+1}d_{m}^{-1}<1\text{,} \\
d_{m+3}d_{m+2}^{-1} &<&1\text{, }d_{m+4}d_{m+3}^{-1}<1\text{,}...\text{,}%
d_{2m+2}d_{2m+1}^{-1}<1\text{,} \\
d_{2m+4}d_{2m+3}^{-1} &<&1\text{, }d_{2m+5}d_{2m+4}^{-1}<1\text{,}...\text{,}%
d_{3m+3}d_{3m+2}^{-1}<1\text{.}
\end{eqnarray*}%
Also,%
\begin{eqnarray*}
c^{-1}d_{1}d_{m+1}^{-1}+c^{-1}d_{1}d_{2m+3}^{-1} &=&c^{-1}\left( \frac{1}{%
1-m\varepsilon }+1\right) \\
&<&c^{-1}\frac{2}{1-m\varepsilon }<1\text{,}
\end{eqnarray*}%
\begin{eqnarray*}
c^{-1}d_{m+2}d_{2m+2}^{-1}+c^{-1}d_{m+2}d_{2m+3}^{-1} &=&c^{-1}\left( \frac{1%
}{1-m\varepsilon }+1\right) \\
&<&c^{-1}\frac{2}{1-m\varepsilon }<1\text{,}
\end{eqnarray*}%
\begin{eqnarray*}
c^{-1}d_{2m+3}d_{m+2}^{-1}+c^{-1}d_{2m+3}d_{3m+3}^{-1} &=&c^{-1}\left( 1+%
\frac{1}{1-m\varepsilon }\right) \\
&<&c^{-1}\frac{2}{1-m\varepsilon }<1\text{.}
\end{eqnarray*}%
Since $B$ has the same eigenvalues as $DBD^{-1}=E=\left( e_{ij}\right) $, we
obtain that%
\begin{eqnarray*}
\max_{1\leq i\leq 3m+3}\left\vert \lambda _{i}\right\vert &\leq &\left\Vert
DBD^{-1}\right\Vert _{\infty } \\
&=&\max_{1\leq i\leq 3m+3}\left\{ \dsum\limits_{j=1}^{3m+3}\left\vert
e_{ij}\right\vert \right\} \\
&=&\max \left\{ 
\begin{array}{c}
d_{2}d_{1}^{-1}\text{, }d_{3}d_{2}^{-1}\text{,}...\text{,}d_{m+1}d_{m}^{-1}%
\text{,} \\ 
d_{m+3}d_{m+2}^{-1}\text{, }d_{m+4}d_{m+3}^{-1}\text{,}...\text{,}%
d_{2m+2}d_{2m+1}^{-1}\text{,} \\ 
d_{2m+4}d_{2m+3}^{-1}\text{, }d_{2m+5}d_{2m+4}^{-1}\text{,}...\text{,}%
d_{3m+3}d_{3m+2}^{-1}\text{,} \\ 
c^{-1}d_{1}d_{m+1}^{-1}+c^{-1}d_{1}d_{2m+3}^{-1}\text{,} \\ 
c^{-1}d_{m+2}d_{2m+2}^{-1}+c^{-1}d_{m+2}d_{2m+3}^{-1}\text{,} \\ 
c^{-1}d_{2m+3}d_{m+2}^{-1}+c^{-1}d_{2m+3}d_{3m+3}^{-1}%
\end{array}%
\right\} \\
&<&1\text{.}
\end{eqnarray*}%
This implies that the equilibrium point of system (\ref{equ:zaq}) is locally
asymptotically stable.

\begin{theorem}
Assume that $A>1$. Then, the positive equilibrium point of system (\ref%
{equ:zaq}) is globally asymptotically stable.
\end{theorem}

\textbf{Proof.} Using Theorem (\ref{teo:jhgfdswedf}), we have%
\begin{eqnarray}
L_{1} &=&\lim_{n\rightarrow \infty }\sup x_{n}\text{, \ }L_{2}=\lim_{n%
\rightarrow \infty }\sup y_{n}\text{, \ }L_{3}=\lim_{n\rightarrow \infty
}\sup z_{n}\text{,}  \label{def:hgfds} \\
m_{1} &=&\lim_{n\rightarrow \infty }\inf x_{n}\text{, }m_{2}=\lim_{n%
\rightarrow \infty }\inf y_{n}\text{, \ }m_{3}=\lim_{n\rightarrow \infty
}\inf z_{n}\text{.}  \notag
\end{eqnarray}%
Then, from (\ref{equ:zaq}) and (\ref{def:hgfds}) we have%
\begin{eqnarray}
L_{1} &\leq &A+\frac{L_{1}}{m_{3}}\text{, \ }L_{2}\leq A+\frac{L_{2}}{m_{3}}%
\text{, \ }L_{3}\leq A+\frac{L_{3}}{m_{2}}\text{,}  \label{ineq:tqwdf} \\
m_{1} &\geq &A+\frac{m_{1}}{L_{3}}\text{, \ }m_{2}\geq A+\frac{m_{2}}{L_{3}}%
\text{, \ }m_{3}\geq A+\frac{m_{3}}{L_{2}}\text{.}  \notag
\end{eqnarray}%
Relations (\ref{ineq:tqwdf}) imply that%
\begin{equation*}
AL_{2}+m_{3}\leq m_{3}L_{2}\leq Am_{3}+L_{2}\text{, \ }AL_{3}+m_{2}\leq
m_{2}L_{3}\leq Am_{2}+L_{3}\text{,}
\end{equation*}%
from which we have%
\begin{equation*}
\left( A-1\right) \left( L_{2}-m_{3}\right) \leq 0\text{, \ }\left(
A-1\right) \left( L_{3}-m_{2}\right) \leq 0\text{.}
\end{equation*}%
Since $A>1$, we get%
\begin{equation*}
L_{2}\leq m_{3}\leq L_{3}\text{, \ }L_{3}\leq m_{2}\leq L_{2}\text{,}
\end{equation*}%
from which%
\begin{equation}
L_{2}=L_{3}=m_{2}=m_{3}\text{.}  \label{equl:kjhg}
\end{equation}%
Moreover, from (\ref{ineq:tqwdf}) it follows that%
\begin{equation*}
L_{1}m_{3}\leq Am_{3}+L_{1}\text{, \ }m_{1}L_{3}\leq AL_{3}+m_{1}\text{,}
\end{equation*}%
from which%
\begin{equation*}
L_{1}\left( m_{3}-1\right) \leq Am_{3}\text{, \ }AL_{3}\leq m_{1}\left(
L_{3}-1\right) \text{.}
\end{equation*}%
Using (\ref{equl:kjhg}), we have%
\begin{equation*}
L_{1}\left( L_{3}-1\right) \leq m_{1}\left( L_{3}-1\right) \text{,}
\end{equation*}%
from which%
\begin{equation*}
L_{1}\leq m_{1}\text{.}
\end{equation*}%
Since $x_{n}$ is bounded, it implys that%
\begin{equation*}
L_{1}=m_{1}\text{.}
\end{equation*}%
Hence, every positive solution $\left\{ \left( x_{n},y_{n},z_{n}\right)
\right\} $ of system (\ref{equ:zaq}) tends to the positive equilibrium
system (\ref{equ:zaq}). So, the proof is completed.

\end{document}